\documentclass[12pt,draft,leqno]{article}
\usepackage{amssymb, eucal, latexsym}

\newtheorem{theorem}{Theorem}[section]
\newtheorem{lm}[theorem]{Lemma}
\newtheorem{exa}[theorem]{Example}

\newtheorem{cor}[theorem]{Corollary}
\newtheorem{pro}[theorem]{Proposition}
\newtheorem{defi}[theorem]{Definition}
\newtheorem{defis}[theorem]{Definitions}
\newtheorem{nota}[theorem]{Notation}

\newtheorem{rem}[theorem]{Remark}

\newtheorem{fact}[theorem]{Fact}

\newtheorem{nist}[theorem]{}

\def\p{\varphi}
\def\a{\alpha}
\def\b{\beta}
\def\d{\delta}

\def\g{\gamma}
\def\GA{\Gamma}

\def\l{\lambda}
\def\LAM{\Lambda}

\def\s{\sigma}

\def\pl{\varphi_\Lambda}

\def\lag{\lambda_A^g}

\def\lra{\longrightarrow}

\def\sbe{\subseteq}
\def\spe{\supseteq}
\def\stm{\setminus}
\def\ems{\emptyset}
\def\nes{\neq\emptyset}

\def\cuk{\,\check{}\,}

\def\ex{\exists}
\def\fa{\forall}

\def\we{\wedge}

\def\bv{\bigvee}

\def\ap{^\prime}
\def\inv{^{-1}}
\def\st{\ |\ }

\def\llx{\ll_{\rho}}
\def\lle{\ll_{\eta}}
\def\lley{\ll_{\eta_Y}}

\def\nin{\not\in}

\def\card #1{\vert #1 \vert}

\def\BB{{\cal B}}
\def\CC{{\cal C}}

\def\KK{{\cal K}}
\def\LL{{\cal L}}

\def\PP{{\cal P}}
\def\TT{{\cal T}}

\def\HLC{{\bf HLC}}
\def\DHLC{{\bf DHLC}}

\def\OHLC{{\bf OHLC}}
\def\DOHLC{{\bf DOHLC}}

\def\2{\mbox{{\bf 2}}}
\def\3{\mbox{{\bf 3}}}

\def\int{\mbox{{\rm int}}}

\def\cl{\mbox{{\rm cl}}}
\def\CL{\mbox{{\rm Clust}}}
\def\BClu{\mbox{{\rm BClust}}}

\def\doc{\hspace{-1cm}{\em Proof.}~~}
\def\sq{\hspace*{\fill} \hbox{\vrule\vbox{\hrule\phantom{o}\hrule}\vrule}}
\def\sqs{\sq \vspace{2mm}}

\def\BBBB{{\rm I}\!{\rm B}}
\def\DDDD{{\rm I}\!{\rm D}}

\title{{\LARGE\bf
Open and other kinds of extensions over local compactifications}\\
\vspace{0.35cm}
{\large\bf Georgi Dimov}\thanks{This paper was supported by the
projects no. 101/2007 $``$Categorical Topology" and no. 005/2009
$``$General and Categorical Topology" of the Sofia
University $``$St. Kl. Ohridski".}\\
\vspace{0.25cm}
 {\footnotesize Dept. of Math. and
Informatics, Sofia University,  Blvd. J. Bourchier 5, 1164 Sofia,
Bulgaria}}

\author{}

\date{}

\begin{document}
\maketitle
\begin{abstract}
{\footnotesize
\noindent

Generalizing de Vries Compactification Theorem \cite{dV2} and
strengthening  Leader Local Compactification Theorem \cite{LE2},
we describe the partially ordered set $(\LL(X),\le)$ of all (up to
equivalence) locally compact Hausdorff extensions of a Tychonoff
space $X$. Using this description, we find the necessary and
sufficient conditions which has to satisfy a map between two
Tychonoff spaces in order to have some kind of extension over
arbitrary given in advance Hausdorff
 local compactifications of these spaces; we
regard the following kinds of extensions: open, quasi-open,
skeletal, perfect, injective, surjective. In this way we
generalize some results of V. Z. Poljakov \cite{Po}.}
\end{abstract}

{\footnotesize {\em  MSC:} primary 54C20, 54D35; secondary 54C10,
54D45, 54E05.

}

\footnotetext[1]{{\footnotesize {\em E-mail address:}
gdimov@fmi.uni-sofia.bg}}

\baselineskip = \normalbaselineskip

\section*{Introduction}

In 1959, V. I. Ponomarev \cite{Pon} proved that if $f:X\lra Y$ is
a perfect open surjection  between two normal Hausdorff spaces $X$
and $Y$ then its extension $\b f:\b X\lra\b Y$ over Stone-\v{C}ech
compactifications of these spaces is an open map; he obtained as
well a more general variant of this theorem which concerns
multi-valued mappings. He posed the following problem:
characterize those continuous maps $f:X\lra Y$ between two
Tychonoff spaces for which the map $\b f$ is open. In 1960, A. D.
Ta\u{i}manov \cite{AT} improved Ponomarev's theorem cited above by
replacing $``$perfect" with $``$closed" (and A. V.
Arhangel'ski\u{i} \cite{AT} generalized Ta\u{i}manov's result for
multi-valued mappings). Later on, V. Z. Poljakov \cite{Po}
described the maps $f:X\lra Y$ between two Tychonoff spaces which
have an {\em open}\/ extension over arbitrary given in advance
Hausdorff compactifications $(cX,c_X)$ and $(cY,c_Y)$ of $X$ and
$Y$ respectively. His work is based on the famous Smirnov
Compactification Theorem \cite{Sm2} which affirms that there
exists an isomorphism between the partially ordered sets (=
posets) $(\KK(X),\le)$ of all, up to equivalence, Hausdorff
compactifications of a Tychonoff space $X$, and $(\PP(X),\preceq)$
of all Efremovi\v{c} proximities on the space $X$; with the help
of this theorem, Ju. M. Smirnov \cite{Sm2} describes  the maps
between two Tychonoff spaces which can be extended {\em
continuously}\/ over arbitrary given in advance compactifications
of these spaces. Analogous assertions for the Hausdorff local
compactifications (= locally compact extensions) of Tychonoff
spaces were proved by S. Leader \cite{LE2}.
%

In this paper we generalize Poljakov's and Leader's theorems and
obtain some other results of this type. We regard the following
kinds of extensions over Hausdorff local compactifications: open,
quasi-open (in the sense of \cite{MP}), perfect, skeletal (in the
sense of \cite{MR}), injective, surjective. We characterize the
functions  between Tychonoff spaces which have  extensions of the
kinds listed above over arbitrary given in advance local
compactifications (see Theorem \ref{zdextcmainl}); in particular,
in Corollary \ref{zdextcmaincbl}, we formulate correctly the
Poljakov's answer to Ponomarev's question (V. Z. Poljakov
\cite{Po} derives it from the general theorem proved by him but
gives in fact only a sufficient condition (see Remark
\ref{rempol})).  The characterizations of all these maps are
obtained here with the help of a strengthening of the Leader Local
Compactification Theorem \cite{LE2} (see Theorems \ref{lead5} and
\ref{zdextcl} below). We give a de Vries-type formulation of the
Leader's theorem (i.e. we describe axiomatically the restrictions
of the Leader's {\em local proximities}\/ on the Boolean algebra
$RC(X)$ of all regular closed subsets of a Tychonoff space $X$)
and prove this new assertion independently of the Leader's theorem
using only our generalization (see \cite{Di3}) of de Vries Duality
Theorem \cite{dV2}. This permits us to use our recent general
results obtained in \cite{Di2,Di5}. Finally, on the base of our
variant of Leader's Theorem,
 we characterize
in the language of {\em local contact algebras}\/ only (i.e.
without mentioning the points of the space)
   the poset $(\LL(X),\le)$ of
all, up to equivalence, Hausdorff local compactifications of $X$,
   where $X$ is  a locally
compact Hausdorff space  (see Theorem \ref{compset}); the algebras
which correspond to the Alexandroff (one-point) compactification
and to the Stone-\v{C}ech compactification of a locally compact
Hausdorff space are described explicitly (see Theorem
\ref{alphabeta}). Let us mention as well that in \cite{DV1} we
described, using the language of non-symmetric proximities, the
surjective continuous maps which have a perfect extension  over
arbitrary given in advance Hausdorff local compactifications.

We now fix the notations.

 If $\CC$ denotes a category, we write
$X\in \card\CC$ if $X$ is
 an object of $\CC$, and $f\in \CC(X,Y)$ if $f$ is a morphism of
 $\CC$ with domain $X$ and codomain $Y$. By $Id_{\CC}$ we denote the identity functor on the category $\CC$.

All lattices are with top (= unit) and bottom (= zero) elements,
denoted respectively by 1 and 0. We do not require the elements
$0$ and $1$ to be distinct.

If $X$ is a set then we denote the power set of $X$ by $P(X)$; the
identity function on $X$ is denoted by $id_X$.

 If
$(X,\tau)$ is a topological space and $M$ is a subset of $X$, we
denote by $\cl_{(X,\tau)}(M)$ (or simply by $\cl(M)$ or
$\cl_X(M)$) the closure of $M$ in $(X,\tau)$ and by
$\int_{(X,\tau)}(M)$ (or briefly by $\int(M)$ or $\int_X(M)$) the
interior of $M$ in $(X,\tau)$.

 If $f:X\lra Y$ is a function and
$M\sbe X$ then $f_{\upharpoonright M}$ is the restriction of $f$
having $M$ as a domain and $f(M)$ as a codomain. Further, we
denote by $\DDDD$  the set of all dyadic numbers of the interval
$(0,1)$ and by $\mathbb{Q}$ the topological space of all rational
numbers with their natural topology.

The  closed maps and the open maps between topological spaces are
assumed to be continuous but are not assumed to be onto. Recall
that a map is {\em perfect}\/ if it is closed and compact (i.e.
point inverses are compact sets).

For all notions and notations not defined here see \cite{Dw, E2,
J, NW}.

\section{Preliminaries}

\begin{defis}\label{conalg}
\rm
An algebraic system $\underline {B}=(B,0,1,\vee,\we, {}^*, C)$ is
called a {\it contact Boolean algebra}\/ or, briefly, {\it contact
algebra} (abbreviated as CA) (\cite{DV2})
 if the system
$(B,0,1,\vee,\we, {}^*)$ is a Boolean algebra (where the operation
$``$complement" is denoted by $``\ {}^*\ $")
  and $C$
is a binary relation on $B$, satisfying the following axioms:

\smallskip

\noindent (C1) If $a\not= 0$ then $aCa$;\\
(C2) If $aCb$ then $a\not=0$ and $b\not=0$;\\
(C3) $aCb$ implies $bCa$;\\
(C4) $aC(b\vee c)$ iff $aCb$ or $aCc$.

\smallskip

\noindent We shall simply write $(B,C)$ for a contact algebra. The
relation $C$  is called a {\em  contact relation}.  If $a\in B$
and $D\sbe B$, we will write $``aCD$" for $``(\fa d\in D)(aCd)$".
We will say that two CA's $(B_1,C_1)$ and $(B_2,C_2)$ are  {\em
CA-isomorphic} iff there exists a Boolean isomorphism $\p:B_1\lra
B_2$ such that, for each $a,b\in B_1$, $aC_1 b$ iff $\p(a)C_2
\p(b)$. A  CA $(B,C)$ is called a {\em complete  contact Boolean
algebra}\/ or, briefly, {\em complete  contact algebra}
(abbreviated as CCA) if $B$ is a complete Boolean algebra.

 A contact algebra $(B,C)$ is called a {\it  normal
contact Boolean algebra}\/ or, briefly, {\it  normal contact
algebra} (abbreviated as NCA) (\cite{dV2,F2}) if it satisfies the
following axioms (we will write $``-C$" for $``not\ C$"):

\smallskip

\noindent (C5) If $a(-C)b$ then $a(-C)c$ and $b(-C)c^*$ for some $c\in B$;\\
(C6) If $a\not= 1$ then there exists $b\not= 0$ such that
$b(-C)a$.

\smallskip

\noindent If an NCA is a CCA then it is called a {\em complete
normal contact Boolean algebra}\/ or, briefly, {\em complete
normal contact algebra} (abbreviated as CNCA).  The notion of
normal contact algebra was introduced by Fedorchuk \cite{F2} under
the name {\em Boolean $\d$-algebra}\/ as an equivalent expression
of the notion of {\em compingent Boolean algebra} of de Vries
\cite{dV2}. We call such algebras $``$normal contact algebras"
because they form a subclass of the class of contact algebras and
naturally arise in the normal Hausdorff spaces.

For any CA $(B,C)$, we define a binary relation  $``\ll_C $"  on
$B$ (called {\em non-tangential inclusion})  by $``\ a \ll_C b
\leftrightarrow a(-C)b^*\ $". Sometimes we will write simply
$``\ll$" instead of $``\ll_C$".
\end{defis}

\begin{exa}\label{extrcr}
\rm Let $B$ be a Boolean algebra. Then there exist a largest and a
smallest contact relations on $B$; the largest one, $\rho_l$, is
defined by $a\rho_l b$ iff $a\neq 0$ and $b\neq 0$, and the
smallest one, $\rho_s$, by $a\rho_s b$ iff $a\wedge b\neq 0$. Note
that, for $a,b\in B$, $a\ll_{\rho_s} b$ iff $a\le b$; hence
$a\ll_{\rho_s} a$, for any $a\in B$. Thus $(B,\rho_s)$ is a normal
contact algebra.
\end{exa}

\begin{exa}\label{rct}
\rm Recall that a subset $F$ of a topological space $(X,\tau)$ is
called {\em regular closed}\/ if $F=\cl(\int (F))$. Clearly, $F$
is regular closed iff it is the closure of an open set. For any
topological space $(X,\tau)$, the collection $RC(X,\tau)$ (we will
often write simply $RC(X)$) of all regular closed subsets of
$(X,\tau)$ becomes a complete Boolean algebra
$(RC(X,\tau),0,1,\we,\vee,{}^*)$ under the following operations: $
1 = X,  0 = \emptyset, F^* = \cl(X\stm F), F\vee G=F\cup G, F\we G
=\cl(\int(F\cap G)). $ The infinite operations are given by the
following formulas: $\bigvee\{F_\g\st
\g\in\GA\}=\cl(\bigcup\{F_\g\st \g\in\GA\}),$ and
$\bigwedge\{F_\g\st \g\in\GA\}=\cl(\int(\bigcap\{F_\g\st
\g\in\GA\})).$

It is easy to see that setting $F \rho_{(X,\tau)} G$ iff $F\cap
G\not = \ems$, we define a contact relation $\rho_{(X,\tau)}$ on
$RC(X,\tau)$; it is called a {\em standard contact relation}. So,
$(RC(X,\tau),\rho_{(X,\tau)})$ is a complete CA (it is called a
{\em standard contact algebra}). We will often write simply
$\rho_X$ instead of $\rho_{(X,\tau)}$. Note that, for $F,G\in
RC(X)$, $F\ll_{\rho_X}G$ iff $F\sbe\int_X(G)$. Clearly, if
$(X,\tau)$ is a normal Hausdorff space then the standard contact
algebra $(RC(X,\tau),\rho_{(X,\tau)})$ is a complete NCA.

A subset $U$ of a topological space $(X,\tau)$ is called {\em
regular open}\/ if $X\stm U\in RC(X)$. The set of all regular open
subsets of $(X,\tau)$ will be denoted by $RO(X,\tau)$ (or simply
by $RO(X)$).
\end{exa}

The next notion and  assertion are inspired by the theory of
proximity spaces (see, e.g., \cite{NW}):

\begin{nist}\label{defcluclan}
\rm Let $(B,C)$ be a CA. Then  a non-empty subset $\s $ of $B$ is
called a {\em cluster in} $(B,C)$
if the following conditions are satisfied:

\smallskip

\noindent (K1) If $a,b\in\s $ then $aCb$;\\
(K2) If $a\vee b\in\s $ then $a\in\s $ or $b\in\s $;\\
(K3) If $aCb$ for every $b\in\s $, then $a\in\s $.

\smallskip

\noindent The set of all clusters in $(B,C)$ will be denoted  by
$\CL(B,C)$.
\end{nist}

\begin{theorem}\label{conclustth}{\rm (\cite{VDDB2})}
A subset $\s$ of a normal contact algebra $(B,C)$ is a cluster iff
there exists an ultrafilter $u$ in $B$ such that
%
$\s=\{a\in B\st aCb \mbox{ for every } b\in u\}.$ Hence, if $u$ is
an ultrafilter  in $B$ then there exists a unique cluster $\s_u$
in $(B,C)$ containing $u$, and
%
$\s_u=\{a\in B\st aCb \mbox{  for every } b\in u\}.$
\end{theorem}

The following notion is a lattice-theoretical counterpart of the
Leader's notion of {\em local proximity} (\cite{LE2}):

\begin{defi}\label{locono}{\rm (\cite{R2})}
\rm An algebraic system $\underline {B}_{\, l}=(B,0,1,\vee,\we,
{}^*, \rho, \BBBB)$ is called a {\it local contact Boolean
algebra}\/ or, briefly, {\it local contact algebra} (abbreviated
as LCA)   if $(B,0,1, \vee,\we, {}^*)$ is a Boolean algebra,
$\rho$ is a binary relation on $B$ such that $(B,\rho)$ is a CA,
and $\BBBB$ is an ideal (possibly non proper) of $B$, satisfying
the following axioms:

\smallskip

\noindent(BC1) If $a\in\BBBB$, $c\in B$ and $a\ll_\rho c$ then
$a\ll_\rho b\ll_\rho c$ for some $b\in\BBBB$;\\
(BC2) If $a\rho b$ then there exists an element $c$ of $\BBBB$
such that
$a\rho (c\we b)$;\\
(BC3) If $a\neq 0$ then there exists  $b\in\BBBB\stm\{0\}$ such
that $b\ll_\rho a$.

\smallskip

We shall simply write  $(B, \rho,\BBBB)$ for a local contact
algebra. When $B$ is a complete Boolean algebra,  the LCA
$(B,\rho,\BBBB)$ is called a {\em complete local contact Boolean
algebra}\/ or, briefly, {\em complete local contact algebra}
(abbreviated as CLCA).

We will say that two local contact algebras $(B,\rho,\BBBB)$ and
$(B_1,\rho_1,\BBBB_1)$ are  {\em LCA-isomorphic} if there exists a
Boolean isomorphism $\p:B\lra B_1$ such that, for $a,b\in B$,
$a\rho b$ iff $\p(a)\rho_1 \p(b)$, and $\p(a)\in\BBBB_1$ iff
$a\in\BBBB$.

Note that if $(B,\rho,\BBBB)$ is a local contact algebra and
$1\in\BBBB$ then $(B,\rho)$ is a normal contact algebra.
Conversely, any normal contact algebra $(B,C)$ can be regarded as
a local contact algebra of the form $(B,C,B)$.
\end{defi}

The following definitions and lemmas
are lattice-theoretical counterparts of some notions and theorems
from Leader's paper \cite{LE2}.

\begin{defi}\label{Alexprn}{\rm (\cite{VDDB2})}
\rm Let $(B,\rho,\BBBB)$ be a local contact algebra. Define a
binary relation $``C_\rho$" on $B$ by
%
$aC_\rho b\ \mbox{ iff }\ a\rho b\ \mbox{ or }\ a,b\not\in\BBBB$;
%
it is called the\/ {\em Alexandroff extension of}\/ $\rho$.
\end{defi}

\begin{lm}\label{Alexprn1}{\rm (\cite{VDDB2})}
Let $(B,\rho,\BBBB)$ be a local contact algebra. Then
$(B,C_\rho)$, where $C_\rho$ is the Alexandroff extension of
$\rho$, is a normal contact algebra.
\end{lm}

\begin{defi}\label{boundcl}
\rm Let $(B,\rho,\BBBB)$ be a local contact algebra. We will say
that $\s$ is a {\em cluster in} $(B,\rho,\BBBB)$ if $\s$ is a
cluster in the NCA $(B,C_\rho)$. A cluster $\s$ in
$(B,\rho,\BBBB)$ (resp., an ultrafilter in $B$) is called {\em
bounded}\/ if $\s\cap\BBBB\nes$ (resp., $u\cap\BBBB\nes$). The set
of all bounded clusters in $(B,\rho,\BBBB)$ will be denoted by
$\BClu(B,\rho,\BBBB)$.
\end{defi}

\begin{lm}\label{neogrn}{\rm (\cite{VDDB2})}
Let $(B,\rho,\BBBB)$ be a local contact algebra and let
$1\not\in\BBBB$. Then $\s_\infty^{(B,\rho,\BBBB)}=\{b\in B\st
b\not\in\BBBB\}$ is a cluster in $(B,\rho,\BBBB)$. (Sometimes we
will simply write $\s_\infty$
instead of $\ \s_\infty^{(B,\rho,\BBBB)}$.)
\end{lm}

\begin{nota}\label{compregn}
\rm Let $(X,\tau)$ be a topological space. We denote by
$CR(X,\tau)$ the family of all compact regular closed subsets of
$(X,\tau)$. We will often write  $CR(X)$ instead of $CR(X,\tau)$.
If $x\in X$ then we set:
\begin{equation}\label{sxvx}
\s_x^X=\{F\in RC(X)\st x\in F\} \mbox{ and } \nu_x^X=\{F\in
RC(X)\st x\in \int_X(F)\}.
\end{equation}
We will often write $\s_x$ and $\nu_x$ instead of, respectively,
$\s_x^X$ and $\nu_x^X$.
\end{nota}

\begin{fact}\label{stanlocn}(\cite{R2,VDDB2})
Let $(X,\tau)$ be a locally compact Hausdorff space. Then:

(a) the triple
$(RC(X,\tau),\rho_{(X,\tau)}, CR(X,\tau))$
  is a complete local contact algebra; it is called a
{\em standard local contact algebra};

(b) for every $x\in X$, $\s_x$ is a bounded cluster in the
standard local contact algebra $(RC(X,\tau),\rho_{(X,\tau)},
CR(X,\tau))$ and $\nu_x$ is a filter in the Boolean algebra
$RC(X)$.
\end{fact}

The next  theorem was proved by Roeper \cite{R2} (but its
particular case concerning compact Hausdorff spaces and NCAs was
proved by de Vries \cite{dV2}).  We will give a sketch of its
proof; it follows the plan of the proof presented in \cite{VDDB2}.
The notations and the facts stated here will be used later on.

\begin{theorem}\label{roeperl}{\rm (P. Roeper \cite{R2})}
 There exists a bijective correspondence between the
class of all (up to isomorphism) CLCAs  and the class of all (up
to homeomorphism) locally compact Hausdorff spaces; its
restriction to the class of all (up to isomorphism) CNCAs gives a
bijective correspondence between the later class and the class of
all (up to homeomorphism) compact Hausdorff spaces.
\end{theorem}

\noindent{\em Sketch of the Proof.}~  Let $(X,\tau)$ be a locally
compact Hausdorff space. We put
%
$$\Psi^t(X,\tau)=(RC(X,\tau),\rho_{(X,\tau)},CR(X,\tau)).$$
%
 Let $(B,\rho,\BBBB)$ be a
complete local contact algebra. Let $C=C_\rho$ be the Alexandroff
extension of $\rho$.  Put $X=\CL(B,C)$ and let $\TT$ be the
topology on $X$ having as a closed base the family
$\{\l_{(B,C)}(a)\st a\in B\}$ where, for every $a\in B$,
%
$$\l_{(B,C)}(a) = \{\s \in X\st  a \in \s\}.$$
%
Sometimes we will write simply $\l_B$ instead of $\l_{(B,C)}$.
%
Note that
%
$X\stm \l_B(a)= \int(\l_B(a^*)),$
%
 the family  $\{\int(\l_B(a))\st a\in B\}$  is an open base of
$(X,\TT)$
%
and, for every $a\in B$,
%
$\l_B(a)\in RC(X,\TT).$
%
%
%
%
Further, $$\l_B:(B,C)\lra (RC(X),\rho_X)$$ is a CA-isomorphism and
%
$(X,\TT)$ is a compact Hausdorff space.
%

Let $1\in\BBBB$. Then $C=\rho$ and $\BBBB=B$, so
that $(B,\rho,\BBBB)=(B,C,B)=(B,C)$ is a
normal contact algebra  and we put
%
$$\Psi^a(B,\rho,\BBBB)(=\Psi^a(B,C,B)=\Psi^a(B,C))=(X,\TT).$$
 Let $1\not\in\BBBB$. Then we set
$L=\BClu(B,\rho,\BBBB)(=X\stm\{\s_\infty\})$ (sometimes we will
write $L_{(B,\rho,\BBBB)}$ or $L_B$ instead of $L$).
 Let the topology $\tau(=\tau_{(B,\rho,\BBBB)})$ on $L$ be the
subspace topology, i.e. $\tau=\TT_{|_L} $. Then $(L,\tau)$ is a
locally compact Hausdorff space. We put
%
$$\Psi^a(B,\rho,\BBBB)=(L,\tau).$$
%
Let
%
$\l^l_{(B,\rho,\BBBB)}(a)=\l_{(B,C_\rho)}(a)\cap L,$
%
for each $a\in B$. We will write simply $\l^l_B$ (or even
$\l_{(B,\rho,\BBBB)}$) instead of $\l^l_{(B,\rho,\BBBB)}$ when
this does not lead to ambiguity. Then
 $L$ is a dense subset of the topological space
$X$ and $\l^l_B: (B,\rho,\BBBB)\lra (RC(L),\rho_L, CR(L))$ is
%
an
LCA-isomorphism.
 Note also that for every $b\in B$,
%
$\int_{L_B}(\l^l_B(b))=L_B\cap\int_X(\l_B(b))$ and
%
%
$L\stm \l_B^l(b)= \int_L(\l_B^l(b^*)).$
%


For every CLCA $(B,\rho,\BBBB)$ and every $a\in B$, set
%
$$\l^g_{(B,\rho,\BBBB)}(a)=\l_{(B,C_\rho)}(a)\cap\Psi^a(B,\rho,\BBBB).$$
%
We will write simply $\l^g_B$ instead of $\l^g_{(B,\rho,\BBBB)}$
when this does not lead to ambiguity.
 Thus, when $1\in\BBBB$,
we have that $\l^g_B=\l_B$, and  if $1\nin\BBBB$  then
$\l^g_B=\l^l_B$. Hence,
%
$$\l^g_B: (B,\rho,\BBBB)\lra (\Psi^t\circ\Psi^a)(B,\rho,\BBBB)$$  is
%
and
%
an LCA-isomorphism.
Let  $(L,\tau)$ be a  locally compact Hausdorff space.  Then
 the map
%
$$t_{(L,\tau)}:(L,\tau)\lra\Psi^a(\Psi^t(L,\tau)),$$
%
defined by  $t_{(L,\tau)}(x)=\s_x$, for all $x\in L$, is a
homeomorphism; we will often write simply $t_L$ instead of
$t_{(L,\tau)}$. Therefore $\Psi^a(\Psi^t(L,\tau))$ is homeomorphic
to $(L,\tau)$ and
$\Psi^t(\Psi^a(B,\rho,\BBBB))$ is LCA-isomorphic to
$(B,\rho,\BBBB)$.
 \sqs

Note that  if $(B,\rho,\BBBB)$ is an LCA, then for every $a\in B$,
%
$a=\bigvee\{b\in\BBBB\st b\ll_{\rho} a\}.$
%

We will need also the following assertion from \cite{Di3}:

\begin{pro}\label{bbcl}{\rm (\cite{Di3})}
Let $(A,\rho,\BBBB)$ be an LCA and $\s_1$, $\s_2$ be two clusters
in $(A,\rho,\BBBB)$ such that\/ $\BBBB\cap\s_1=\BBBB\cap\s_2$.
Then $\s_1=\s_2$.
\end{pro}


We will recall some results from \cite{Di3,Di2} which are basic
for our investigations in the present paper.

\begin{defis}\label{dvfi}{\rm (\cite{Di3})}
\rm Let  $\HLC$ be the category of all locally compact Hausdorff
spaces and all continuous maps between them.

Let $\DHLC$ be the category whose objects are all complete LCAs
and whose morphisms are all functions $\p:(A,\rho,\BBBB)\lra
(B,\eta,\BBBB\ap)$ between the objects of $\DHLC$ satisfying
conditions

\smallskip

\noindent(DLC1) $\p(0)=0$;\\
(DLC2) $\p(a\we b)=\p(a)\we \p(b)$, for all $a,b\in A$;\\
(DLC3) If $a\in\BBBB, b\in A$ and $a\llx b$, then $(\p(a^*))^*\lle
\p(b)$;\\
(DLC4) For every $b\in\BBBB\ap$ there exists $a\in\BBBB$ such that
$b\le\p(a)$;\\
\noindent(DLC5) $\p(a)=\bigvee\{\p(b)\st b\in\BBBB, b\llx a\}$,
for every $a\in A$;

\medskip

{\noindent}let the composition $``\diamond$" of two morphisms
$\p_1:(A_1,\rho_1,\BBBB_1)\lra (A_2,\rho_2,\BBBB_2)$ and
$\p_2:(A_2,\rho_2,\BBBB_2)\lra (A_3,\rho_3,\BBBB_3)$ of\/ $\DHLC$
be defined by the formula
%
$\p_2\diamond\p_1 = (\p_2\circ\p_1)\cuk,$
%
 where, for every
function $\psi:(A,\rho,\BBBB)\lra (B,\eta,\BBBB\ap)$ between two
objects of\/ $\DHLC$, $\psi\cuk:(A,\rho,\BBBB)\lra
(B,\eta,\BBBB\ap)$ is defined as follows:
%
$\psi\cuk(a)=\bigvee\{\psi(b)\st b\in \BBBB, b\llx a\},$
%
for every $a\in A$.
\end{defis}

 As it was shown in \cite{Di3},  condition (DLC3) can be
replaced by the following one:

\smallskip

\noindent(DLC3S) If $a, b\in A$ and $a\llx b$, then
$(\p(a^*))^*\lle \p(b)$.

\smallskip

We will need the following duality theorem:

\begin{theorem}\label{lccont}{\rm (\cite{Di3})}
The categories $\HLC$ and\/ $\DHLC$ are dually equivalent. In more
details, let
$\LAM^t:\HLC\lra\DHLC$ and $\LAM^a:\DHLC\lra\HLC$
be the contravariant functors extending, respectively,  the
Roeper's correspondences $\Psi^t$  and\/ $\Psi^a$ (see Theorem
\ref{roeperl}) to the
 morphisms of the categories $\HLC$ and\/ $\DHLC$ in
the following way: for every $f\in\HLC(X,Y)$ and every $G\in
RC(Y)$,
$$\LAM^t(f)(G)=\cl(f\inv(\int(G))),$$ and for every
$\p\in\DHLC((A,\rho,\BBBB),(B,\eta,\BBBB\ap))$ and for every
$\s\ap\in\LAM^a(B,\eta,\BBBB\ap)$,
\begin{equation}\label{lamdefs}
\LAM^a(\p)(\s\ap)\cap\BBBB=\{a\in \BBBB\st \mbox{if } b\in A
\mbox{ and } a\llx b \mbox{ then }\p(b)\in\s\ap\}
\end{equation}
(if, in addition, $\p$ is a complete Boolean homomorphism, then
the above formula is equivalent to the following one: for every
bounded ultrafilter $u$ in $B$,
$\LAM^a(\p)(\s_u)=\s_{\p\inv(u)}$); then
$\l^g: Id_{\,\DHLC}\lra\LAM^t\circ\LAM^a,$ where
$\l^g(A,\rho,\BBBB)=\l_A^g$
for every $(A,\rho,\BBBB)\in\card\DHLC$ , and
$t^l:Id_{\,\HLC}\lra\LAM^a\circ\LAM^t,$ where $t^l(X)=t_X$
for every $X\in\card\HLC$, are natural isomorphisms.
\end{theorem}

\begin{defi}\label{lcat}{\rm (\cite{Di2})}
\rm Let $\OHLC$ be the category of all locally compact Hausdorff
spaces and all open maps between them.

Let $\DOHLC$ be the subcategory of the category $\DHLC$ having the
same objects and whose morphisms are all $\DHLC$-morphisms
$\p:(A,\rho,\BBBB)\lra (B,\eta,\BBBB\ap)$ which are complete
Boolean homomorphisms and satisfy the following condition:

\smallskip

\noindent (LO) $\fa a\in A$ and $\fa b\in \BBBB\ap$,
$\p_\LAM(b)\rho a$ implies
 $b\eta\p(a)$,

\smallskip

\noindent where $\pl$ is the left adjoint of $\p$ (i.e. $\pl:B\lra
A$ is an order-preserving map such that $\fa b\in B$, $\p(\p_\LAM
(b))\ge b$ and
 $\fa a\in A$, $\p_\LAM (\p(a))\le a$; its existence follows from the Adjoint
Functor Theorem (see, e.g., \cite{J})).
\end{defi}

\begin{theorem}\label{maintheorem}{\rm (\cite{Di2})}
The categories $\OHLC$ and $\DOHLC$ are dually equivalent.
\end{theorem}


Finally, we will recall some definitions and facts from the theory
of extensions of topological spaces, as well as the fundamental
Leader Local Compactification Theorem \cite{LE2}.

 Let $X$ be a Tychonoff space. We will denote by $\LL(X)$ the
set of all, up to equivalence, locally compact Hausdorff
extensions of $X$ (recall that two (locally compact Hausdorff)
extensions $(Y_1,f_1)$ and $(Y_2,f_2)$ of $X$ are said to be {\em
equivalent}\/ iff there  exists a homeomorphism $h:Y_1\lra Y_2$
such that $h\circ f_1=f_2$); the equivalence class of $(Y,f)$ will
be denoted by $[(Y,f)]$.  Let $[(Y_i,f_i)]\in\LL(X)$, where
$i=1,2$. We set $[(Y_1,f_1)]\le [(Y_2,f_2)]$
 if there exists a
continuous  mapping $h:Y_2\lra Y_1$ such that $f_1=h\circ f_2$.
Then $(\LL(X),\le)$ is a poset.

 Let $X$ be a Tychonoff space. We will denote by $\KK(X)$ the
set of all, up to equivalence,  Hausdorff compactifications of
$X$.

Recall that if $X$ is a set and $P(X)$ is the power set of $X$
ordered by the inclusion (and thus $P(X)$ becomes a Boolean
algebra), then a triple $(X,\b,\BB)$ is called a {\em local
proximity space} (see \cite{LE2}) if $(P(X),\b)$ is a CA, $\BB$ is
an ideal (possibly non proper) of $P(X)$ and the axioms
(BC1),(BC2) from \ref{locono}
 are fulfilled. A local proximity
space $(X,\b,\BB)$ is said to be {\em separated} if $\b$ is the
identity relation on singletons. Recall that every separated local
proximity space $(X,\b,\BB)$ induces a Tychonoff topology
$\tau_{(X,\b,\BB)}$ in $X$ by defining $\cl(M)=\{x\in X\st x\b
M\}$ for every $M\sbe X$ (\cite{LE2}). If $(X,\tau)$ is a
topological space then we say that $(X,\b,\BB)$ is a {\em local
proximity space on} $(X,\tau)$ if $\tau_{(X,\b,\BB)}=\tau$.

The  set of all separated local proximity spaces on a Tychonoff
space $(X,\tau)$ will be denoted by $\LL\PP(X,\tau)$. A partial
order in $\LL\PP(X,\tau)$ is defined by $(X,\b_1,\BB_1)\preceq
(X,\b_2,\BB_2)$ if $\b_2\sbe\b_1$ and $\BB_2\sbe\BB_1$ (see
\cite{LE2}).

A function $f:X_1\lra X_2$ between two local proximity spaces
$(X_1,\b_1,\BB_1)$ and $(X_2,\b_2,\BB_2)$ is said to be an  {\em
equicontinuous mapping\/} (see \cite{LE2}) if the following two
conditions are fulfilled:
\smallskip

\noindent(EQ1) $A\b_1 B$ implies $f(A)\b_2 f(B)$, for $A,B\sbe X$,
and
\smallskip

\noindent(EQ2) $B\in\BB_1$ implies $f(B)\in\BB_2$.

The separated local proximity spaces of the form $(X,\d,P(X))$ are
denoted by $(X,\d)$ and are called {\em Efremovi\v{c} proximity
spaces}. The equicontinuous mappings between Efremovi\v{c}
proximity spaces are called {\em proximally continuous mappings}.

\begin{theorem}\label{Leader} {\rm (S. Leader \cite{LE2})}
Let $(X,\tau)$ be a Tychonoff space. Then there exists an
isomorphism $\LAM_X$ between the ordered sets $(\LL(X,\tau),\le)$
and $(\LL\PP(X,\tau),\preceq)$.  In more details, for every $(X,
\b, \BB)\in\LL\PP(X,\tau)$ there exists a locally compact
Hausdorff extension $(Y,f)$  of X satisfying the following two
conditions:

\smallskip

\noindent(a) $A \b B$ iff $\cl_Y(f(A))\cap \cl_Y(f(B))\nes$;

\smallskip

\noindent(b) $B\in\BB$  iff $\cl_Y(f(B))$ is compact.

\smallskip

\noindent Such a local compactification is unique up to
equivalence; we set $(Y,f)=L(X,\b,\BB)$ and
$(\LAM_X)\inv(X,\b,\BB)=[(Y,f)]$. The space $Y$ is compact iff
$X\in\BB$. Conversely, if $(Y,f)$ is a locally compact Hausdorff
extension of $X$ and $\b$ and $\BB$ are defined by (a) and (b),
then $(X, \b, \BB)$ is a separated local proximity space, and we
set $\LAM_X([(Y,f)])=(X,\b,\BB)$.

Let $(X_i,\b_i,\BB_i)$, $i=1,2$, be two separated local proximity
spaces and $f:X_1\lra X_2$ be a function. Let
$(Y_i,f_i)=L(X_i,\b_i,\BB_i)$, where $i=1,2$. Then there exists a
continuous map $L(f):Y_1\lra Y_2$ such that $f_2\circ f= L(f)\circ
f_1$ iff $f$ is an equicontinuous map between $(X_1,\b_1,\BB_1)$
and $(X_2,\b_2,\BB_2)$.
\end{theorem}

 We will also need a lemma from \cite{CNG}:

\begin{lm}\label{isombool}
Let $X$ be a dense subspace of a topological space $Y$. Then the
functions $r:RC(Y)\lra RC(X)$, $F\mapsto F\cap X$, and
$e:RC(X)\lra RC(Y)$, $G\mapsto \cl_Y(G)$, are Boolean isomorphisms
between Boolean algebras $RC(X)$ and $RC(Y)$, and $e\circ
r=id_{RC(Y)}$, $r\circ e=id_{RC(X)}$.
\end{lm}


\section{A de Vries-type revision of the Leader Local Compactification Theorem}

In this section we will obtain a strengthening of  Leader Local
Compactification
 Theorem  (\cite{LE2}); it is similar to de Vries'
(\cite{dV2})
 strengthening of Smirnov Compactification Theorem (\cite{Sm2}).

\begin{defi}\label{lead2}
\rm  Let $(X,\tau)$ be a Tychonoff space. An LCA
$(RC(X,\tau),\rho,\BBBB)$ is said to be {\em admissible for}
$(X,\tau)$ if it satisfies the following conditions:

\smallskip

\noindent(A1) if $F,G\in RC(X)$ and $F\cap G\nes$ then $F\rho G$;

\smallskip

\noindent(A2) if $F\in RC(X)$ and $x\in\int_X(F)$ then there
exists $G\in\BBBB$ such that $x\in\int_X(G)$ and $G\llx F$.

\smallskip

The set of all LCAs $(RC(X,\tau),\rho,\BBBB)$  which are
admissible for $(X,\tau)$ will be denoted by $\LL_{ad}(X,\tau)$
(or simply by $\LL_{ad}(X)$). If
$(RC(X),\rho_i,\BBBB_i)\in\LL_{ad}(X)$, where $i=1,2$, then we set
$(RC(X),\rho_1,\BBBB_1)\preceq_{ad} (RC(X),\rho_2,\BBBB_2)$ iff
$\rho_2\sbe\rho_1$ and $\BBBB_2\sbe\BBBB_1$. Obviously,
$(\LL_{ad}(X,\tau),\preceq_{ad})$ is a poset.
\end{defi}

\begin{theorem}\label{lead5}
Let $(X,\tau)$ be a Tychonoff space. Then the
posets $(\LL(X,\tau),\le)$ and $(\LL_{ad}(X,\tau),\preceq_{ad})$
are isomorphic.
\end{theorem}

\doc
Let $(Y,f)$ be a locally compact Hausdorff extensions of $X$. Set
\begin{equation}\label{01l}
\BBBB_{(Y,f)}=f\inv(CR(Y)) \mbox{  and let } F\eta_{(Y,f)} G\iff
\cl_Y(f(F))\cap\cl_Y(f(G))\nes,
\end{equation}
 for every $F,G\in RC(X)$. Note
that, by \ref{isombool}, $\BBBB_{(Y,f)}=\{F\in RC(X)\st
\cl_Y(f(F))$ is compact$\}$. Hence $\BBBB_{(Y,f)}\sbe RC(X)$. We
will show that $(RC(X),\eta_{(Y,f)},\BBBB_{(Y,f)})\in\LL_{ad}(X)$.
We have, by \ref{isombool}, that the map
\begin{equation}\label{001l}
r_{(Y,f)}:(RC(Y),\rho_Y,CR(Y))\lra(RC(X),\eta_{(Y,f)},\BBBB_{(Y,f)}),
\  G\mapsto f\inv(G),
\end{equation}
is a Boolean isomorphism  and, for every $F,G\in RC(Y)$, the
following is fulfilled: $F\rho_{(Y,f)} G$ iff
$r_{(Y,f)}(F)\eta_{(Y,f)} r_{(Y,f)}(G)$, and  $F\in CR(Y)$ iff
$r_{(Y,f)}(F)\in\BBBB_{(Y,f)}$. Hence
$(RC(X),\eta_{(Y,f)},\BBBB_{(Y,f)})$ is an LCA and $r_{(Y,f)}$ is
an LCA-isomorphism. Clearly, condition (A1) is fulfilled. Let now
$F\in RC(X)$. Set $U=\int_X(F)$ and let $x\in U$. There exists an
open subset $V$ of $Y$ such that $V\cap f(X)=f(U)$. Since $Y$ is a
locally compact Hausdorff space, there exists an $H\in CR(Y)$ with
$f(x)\in\int_Y(H)\sbe H\sbe V$. Let $G=f\inv(H)$. Then
$H\in\BBBB_{(Y,f)}$ and, obviously, $x\in\int_X(G) $ and $G\lley
F$. So, condition (A2) is also checked. Hence
$(RC(X),\eta_{(Y,f)},\BBBB_{(Y,f)})\in\LL_{ad}(X)$. It is clear
that  if $(Y_1,f_1)$ is a locally compact Hausdorff extensions of
$X$ equivalent to the extension $(Y,f)$, then
$(RC(X),\eta_{(Y,f)},\BBBB_{(Y,f)})=(RC(X),\eta_{(Y_1,f_1)},\BBBB_{(Y_1,f_1)})$.
Therefore, a map
\begin{equation}\label{dw1l}
\a_X:\LL(X)\lra\LL_{ad}(X), \
[(Y,f)]\mapsto(RC(X),\eta_{(Y,f)},\BBBB_{(Y,f)}),
\end{equation}
 is well-defined.

Set, for short, $A=RC(X)$. Let $(A,\rho,\BBBB)\in\LL_{ad}(X)$ and
$Y=\LAM^a(A,\rho,\BBBB)$. Then, by Roeper's Theorem \ref{roeperl},
$Y$ is a locally compact Hausdorff space. Let us show that for
every $x\in X$, we have that $\s_x\in Y$ (where $\s_x=\{F\in A\st
x\in F\}$).  By \ref{stanlocn}, $\nu_x$ is a filter in the Boolean
algebra $A$. Hence there exists an ultrafilter $u$ in $A$ such
that $\nu_x\sbe u$. It is easy to see that $u\sbe\s_x$. Let
$\s=\{F\in A\st FC_\rho u\}$ (i.e. $\s=\s_u$). Since, by (A2),
$\nu_x\cap\BBBB\nes$, we get that $\s\in Y$. We will show that
$\s_x=\s$. Indeed, let $F\in\s_x$ and $G\in u$. Then $x\in F\cap
G$. Thus, by (A1), $F\rho G$. This implies that $FC_\rho u$, i.e.
that $F\in\s$. So, $\s_x\sbe\s$. Now, suppose that there exists
$F\in\s$ such that $x\nin F$. Then $x\in X\stm F=\int_X(F^*)$.
Thus, by (A2), there exists $G\in\BBBB$ such that $x\in\int_X(G)$
and $G\llx F^*$. Therefore $G\in\nu_x$ and $G(-\rho)F$. Since
$G\in\BBBB$, we get that $F(-C_\rho)G$, a contradiction. So, we
have proved that $\s_x=\s$ and, thus, $\s_x\in Y$ for every $x\in
X$. Define
\begin{equation}\label{0fl}
f_{(\rho,\BBBB)}:X\lra Y, \  x\mapsto\s_x.
\end{equation}
 Set, for short,
$f=f_{(\rho,\BBBB)}$. Then $\cl_Y(f(X))=Y$. Indeed, for every
$F\in\BBBB\stm\{\ems\}$ and for every $x\in F$, we have that
$\s_x\in f(X)\cap\lag(F)$. Since $Y$ is regular, this implies that
$\cl_Y(f(X))=Y$. We will now show that $f$ is a homeomorphic
embedding. It is clear that $f$ is an injection. Further, let
$x\in X$, $F\in\BBBB$ and $\s_x\in\int_Y(\lag(F))$. Since
$\int_Y(\lag(F))=Y\stm\lag(F^*)$, we get that $\s_x\nin\lag(F^*)$.
Thus $F^*\nin\s_x$. This implies that $x\nin F^*$, i.e. $x\in
X\stm F^*=\int_X(F)$. Moreover, $f(\int_X(F))\sbe\int_Y(\lag(F))$.
Indeed, if $y\in\int_X(F)$ then $y\nin F^*$; thus $F^*\nin\s_y$,
i.e. $\s_y\nin\lag(F^*)$; this implies that
$\s_y\in\int_Y(\lag(F))$. All this shows that $f$ is a continuous
function. Set $g=((f)_{\upharpoonright X})\inv$, where
$(f)_{\upharpoonright X}:X\lra f(X)$ is the restriction of $f$. We
have to show that $g$ is a continuous function. Let $x\in X$,
$F\in A$ and $x\in\int_X(F)$. We have that $x=g(\s_x)$. Let
$\s_y\in\int_Y(\lag(F))$. Then $\s_y\in Y\stm\lag(F^*)$, i.e.
$y\nin F^*$; thus $y\in X\stm F^*=\int_X(F)$. Therefore,
$g(\int_Y(\lag(F)))\sbe\int_X(F)$. So, $g$ is a continuous
function. All this shows that $(Y,f)$ is a locally compact
Hausdorff extension of $X$. We now set:
\begin{equation}\label{dw2l}
\b_X:\LL_{ad}(X)\lra \LL(X), \  (RC(X),\rho,\BBBB)\mapsto
[(\LAM^a(RC(X),\rho,\BBBB),f_{(\rho,\BBBB)})].
\end{equation}

We will show that $\a_X\circ\b_X=id_{\LL_{ad}(X)}$ and
$\b_X\circ\a_X=id_{\LL(X)}$.

Let $[(Y,f)]\in\LL(X)$. Then
$\b_X(\a_X([(Y,f)]))=\b_X(RC(X),\eta_{(Y,f)},\BBBB_{(Y,f)})=
[(\LAM^a(RC(X),\eta_{(Y,f)},\BBBB_{(Y,f)}),f_{(\eta_{(Y,f)},\BBBB_{(Y,f)})})]$.
Set, for short, $\eta=\eta_{(Y,f)}$, $\BBBB=\BBBB_{(Y,f)}$,
$g=f_{(\eta_{(Y,f)},\BBBB_{(Y,f)})}$,
$Z=\LAM^a(RC(X),\eta_{(Y,f)},\BBBB_{(Y,f)})$ and $r_{(Y,f)}=r_f$.
We have to show that $[(Y,f)]=[(Z,g)]$. Since $r_f$ is an
LCA-isomorphism, we get that
$h=\LAM^a(r_f):Z\lra\LAM^a(\LAM^t(Y))$ is a homeomorphism. Set
$Y\ap=\LAM^a(\LAM^t(Y))$.  By Roeper's Theorem \ref{roeperl}, the
map $t_Y:Y\lra Y\ap, \  y\mapsto\s_y$ is a homeomorphism. Let
$h\ap=(t_Y)\inv\circ h$. Then $h\ap:Z\lra Y$ is a homeomorphism.
We will prove that $h\ap\circ g=f$ and this will imply that
$[(Y,f)]=[(Z,g)]$. Let $x\in X$ and $u$ be an ultrafilter
containing the filter $\nu_x$. Then, as we have shown above,
$\s_x=\s_u$. Hence $h(\s_x)=h(\s_u)=\s_{e_f(u)}$, where
$e_f=(r_f)\inv$. Thus
$h\ap(g(x))=h\ap(\s_x)=(t_Y)\inv(h(\s_x))=(t_Y)\inv(\s_{e_f(u)})$.
Note that, by \ref{isombool}, $e_f(F)=\cl_Y(f(F))$, for every
$F\in RC(X)$. Since $e_f:RC(X)\lra RC(Y)$ is a Boolean
isomorphism, we get that $e_f(u)$ is an ultrafilter in $RC(Y)$
containing $\nu^Y_{f(x)}$. Thus $\s_{e_f(u)}=\s^Y_{f(x)}$. Hence
$(t_Y)\inv(\s_{e_f(u)})=f(x)$. So, $h\ap\circ g=f$. Therefore,
$\b_X\circ\a_X=id_{\LL(X)}$.

Let $(RC(X),\rho,\BBBB)\in\LL_{ad}(X)$ and
$Y=\LAM^a(RC(X),\rho,\BBBB)$. Recall that we have set  $A=RC(X)$.
We have that $\b_X(A,\rho,\BBBB)=[(Y,f_{(\rho,\BBBB)})]$.
 Set $f=f_{(\rho,\BBBB)}$. Then $\a_X(\b_X(A,\rho,\BBBB))=(A,\eta_{(Y,f)},\BBBB_{(Y,f)})$.
   By Roeper's Theorem \ref{roeperl}, we have that
$\lag:(A,\rho,\BBBB)\lra(RC(Y),\rho_Y,CR(Y))$ is an
LCA-isomorphism. We will show that $f\inv(\lag(F))=F$, for every
$F\in RC(X)$. Indeed, if $x\in F$ then $F\in\s_x$, and thus
$\s_x\in\lag(F)$; hence $f(F)\sbe\lag(F)$, i.e. $F\sbe
f\inv(\lag(F))$. If $x\in f\inv(\lag(F))$ then $f(x)\in\lag(F)$,
i.e. $\s_x\in\lag(F)$; therefore $F\in\s_x$, which means that
$x\in F$. So, $f\inv(\lag(F))=F$, for every $F\in RC(X)$. Since
$CR(Y)=\{\lag(F)\st F\in\BBBB\}$, we get that
$f\inv(CR(Y))=\BBBB$. Thus $\BBBB_{(Y,f)}=\BBBB$. Further, by
\ref{isombool}, $\cl_Y(f(F))=\lag(F)$, for every $F\in RC(X)$.
Since, for every $F,G\in RC(X)$, $F\rho G\iff
\lag(F)\cap\lag(G)\nes$, we get that $\rho=\eta_{(Y,f)}$.
Therefore, $\a_X\circ\b_X=id_{\LL_{ad}(X)}$.

We will now prove that $\a_X$ and $\b_X$ are monotone maps.

Let $[(Y_i,f_i)]\in\LL(X)$, where $1=1,2$, and
$[(Y_1,f_1)]\le[(Y_2,f_2)]$. Then there exists a continuous map
$g:Y_2\lra Y_1$ such that $g\circ f_2=f_1$. Let
$\a_X([(Y_i,f_i)])=(RC(X),\eta_{(Y_i,f_i)},\BBBB_{(Y_i,f_i)})$,
where $i=1,2$. Set $\eta_i=\eta_{(Y_i,f_i)}$ and
$\BBBB_i=\BBBB_{(Y_i,f_i)}$,  $i=1,2$. We have to show that
$\eta_2\sbe\eta_1$ and $\BBBB_2\sbe\BBBB_1$. Let $F\in\BBBB_2$.
Then $\cl_{Y_2}(f_2(F))$ is compact. Hence $g(\cl_{Y_2}(f_2(F)))$
is compact. We have that $f_1(F)=g(f_2(F))\sbe
g(\cl_{Y_2}(f_2(F)))\sbe\cl_{Y_1}(g(f_2(F)))=\cl_{Y_1}(f_1(F))$.
Thus $\cl_{Y_1}(f_1(F))=g(\cl_{Y_2}(f_2(F)))$, i.e.
$\cl_{Y_1}(f_1(F))$ is compact. Therefore $F\in\BBBB_1$. So, we
have proved that $\BBBB_2\sbe\BBBB_1$. Let $F,G\in RC(X)$ and
$F\eta_2 G$. Then there exists $y\in
\cl_{Y_2}(f_2(F))\cap\cl_{Y_2}(f_2(G))$. Since
$g(\cl_{Y_2}(f_2(F)))\sbe\cl_{Y_1}(f_1(F))$ and, analogously,
$g(\cl_{Y_2}(f_2(G)))\sbe\cl_{Y_1}(f_1(G))$, we get that
$g(y)\in\cl_{Y_1}(f_1(F))\cap\cl_{Y_1}(f_1(G))$. Thus $F\eta_1 G$.
Therefore, $\eta_2\sbe\eta_1$. All this shows that
$\a_X([(Y_1,f_1)])\preceq_{ad}\a_X([(Y_2,f_2)])$. Hence, $\a_X$ is
a monotone function.

Let now $(RC(X),\rho_i,\BBBB_i)\in\LL_{ad}(X)$, where $i=1,2$, and
$(RC(X),\rho_1,\BBBB_1)\preceq_{ad}(RC(X),\rho_2,\BBBB_2)$. Set,
for short, $Y_i=\LAM^a(RC(X),\rho_i,\BBBB_i)$ and
$f_i=f_{(\rho_i,\BBBB_i)}$, $i=1,2$. Then
$\b_X(RC(X),\rho_i,\BBBB_i)=[(Y_i,f_i)]$, $i=1,2$. We will show
that $[(Y_1,f_1)]\le[(Y_2,f_2)]$. We have that $f_i:X\lra Y_i$ is
defined by $f_i(x)=\s_x$, for every $x\in X$ and $i=1,2$. We also
have that $\BBBB_2\sbe\BBBB_1$ and $\rho_2\sbe\rho_1$. Let us
regard the following function
$\p:(RC(X),\rho_1,\BBBB_1)\lra(RC(X),\rho_2,\BBBB_2), \  F\mapsto
F.$
We will prove that $\p$ is a $\DHLC$-morphism. Clearly, $\p$
satisfies conditions (DLC1) and (DLC2). The fact that
$\rho_2\sbe\rho_1$ implies immediately that $\p$ satisfies also
condition (DLC3). Further, for establishing condition (DLC4) use
the fact that $\BBBB_2\sbe\BBBB_1$. Let $F\in RC(X)$. Then
$F=\bv\{G\in\BBBB_1\st G\ll_{\rho_1}F\}$ and thus
$\p(F)=\bv\{\p(G)\st G\in\BBBB_1, G\ll_{\rho_1}F\}$.  This shows
that $\p$ satisfies condition (DLC5). So, $\p$ is a
$\DHLC$-morphism. Then $g=\LAM^a(\p):Y_2\lra Y_1$ is a continuous
map. We will prove that $g\circ f_2=f_1$, i.e. that for every
$x\in X$, $g(\s_x)=\s_x$. So, let $x\in X$. We have, by
(\ref{lamdefs}),
 that
$g(\s_x)\cap\BBBB_1=\{F\in\BBBB_1\st (\fa G\in
RC(X))[(F\ll_{\rho_1}G)\rightarrow(x\in G)]\}$. We will show that
$g(\s_x)\cap\BBBB_1=\s_x\cap\BBBB_1$. This will imply, by
\ref{bbcl}, that $g(\s_x)=\s_x$. Let $F\in\s_x\cap \BBBB_1$. Then
$x\in F$ and thus  $F\in g(\s_x)\cap\BBBB_1$. Conversely, suppose
that there exists $H\in g(\s_x)\cap\BBBB_1$ such that $x\nin H$.
Then $x\in X\stm H=\int_X(H^*)$. By (A2), there exists
$G\in\BBBB_1$ with $x\in\int_X(G)$ and $G\ll_{\rho_1}H^*$. We get
that $H\ll_{\rho_1}G^*$ and $x\nin G^*$, a contradiction.
Therefore, $g(\s_x)=\s_x$. Thus $[(Y_1,f_1)]\le[(Y_2,f_2)]$. So,
$\b_X$ is also a monotone function. Since $\b_X=(\a_X)\inv$, we
get that $\a_X$  is an isomorphism.
 \sqs

\begin{defi}\label{lead6}
\rm Let $(X,\tau)$ be a Tychonoff space. An NCA $(RC(X,\tau),C)$
is said to be {\em admissible for} $(X,\tau)$ if the LCA
$(RC(X,\tau),C,RC(X,\tau))\in\LL_{ad}(X,\tau)$. The set of all
NCAs which are admissible for $(X,\tau)$ will be denoted by
$\KK_{ad}(X,\tau)$ (or simply by $\KK_{ad}(X)$). Note that
$\KK_{ad}(X)$ is, in fact, a subset of $\LL_{ad}(X)$. The
restriction on $\KK_{ad}(X)$ of the order $\preceq_{ad}$, defined
on $\LL_{ad}(X)$, will be denoted again by $\preceq_{ad}$.
\end{defi}

\begin{cor}\label{lead7}{\rm (de Vries \cite{dV2})}
For every Tychonoff space $X$,
there exists an isomorphism between the posets $(\KK(X),\le)$ and
$(\KK_{ad}(X),\preceq_{ad})$.
\end{cor}

\doc It follows immediately from Theorem \ref{lead5}. \sqs

The first part of Leader Local Compactification Theorem
\ref{Leader} follows from our Theorem \ref{lead5} and the
following three lemmas.

\begin{lm}\label{lead1}
Let $(X,\b_i,\BB_i)$, $i=1,2$, be two separated local proximity
spaces on a Tychonoff space $(X,\tau)$,  $\BB_1\cap
RC(X)=\BB_2\cap RC(X)$ and $(\b_1)_{|RC(X)}= (\b_2)_{|RC(X)}$
(i.e., for every $F,G\in RC(X)$, $F\b_1 G\iff F\b_2 G$). Then
$\b_1=\b_2$ and $\BB_1=\BB_2$.
\end{lm}

\doc Let, for  $i=1,2$, $(Y_i,f_i)=L(X,\b_i,\BB_i)$   (see Theorem \ref{Leader}). Let
$B\in\BB_1$. Then $\cl_{Y_1}(f_1(B))$ is compact. There exists an
open subset $U$ of $Y_1$ such that $\cl_{Y_1}(f_1(B))\sbe U$ and
$\cl_{Y_1}(U)$ is compact. Let $F=f_1\inv(\cl_{Y_1}(U))$. Then
$F\in RC(X)$ and $\cl_{Y_1}(f_1(F))=\cl_{Y_1}(U)$. Hence
$F\in\BB_1\cap RC(X)$. Thus $F\in\BB_2$. Since $B\sbe F$, we get
that $B\in\BB_2$. Therefore, $\BB_1\sbe\BB_2$. Analogously we
obtain that $\BB_2\sbe \BB_1$. Thus $\BB_1=\BB_2$.

Let $M,N\sbe X$ and $M(-\b_1)N$. Suppose that $M\b_2 N$. Then
there exist $M\ap,N\ap\in\BB_2$ such that $M\ap\sbe M$, $N\ap\sbe
N$ and $M\ap\b_2 N\ap$. Since $\BB_1=\BB_2$, we get that
$M\ap,N\ap\in\BB_1$. Hence $K_1=\cl_{Y_1}(f_1(M\ap))$ and
$K_2=\cl_{Y_1}(f_1(N\ap))$ are disjoint compact subsets of $Y_1$.
Then there exist open subsets $U$ and $V$ of $Y_1$ having disjoint
closures in $Y_1$ and containing, respectively, $K_1$ and $K_2$.
Set $F=f_1\inv(\cl_{Y_1}(U))$ and $G=f_1\inv(\cl_{Y_1}(V))$. Then
$F,G\in RC(X)$, $M\ap\sbe F$, $N\ap\sbe G$ and $F(-\b_1)G$. Thus
$F(-\b_2)G$ and hence $M\ap(-\b_2)N\ap$, a contradiction.
Therefore, $M(-\b_2)N$. So, $\b_2\sbe\b_1$. Using the symmetry, we
obtain that $\b_1=\b_2$. \sqs

\begin{lm}\label{lead3}
Let $(X,\b,\BB)$ be a separated local proximity space. Set
$\tau=\tau_{(X,\b,\BB)}$. Let $\rho=\b_{|RC(X,\tau)}$ and
$\BBBB=\BB\cap RC(X,\tau)$. Then
$(RC(X,\tau),\rho,\BBBB)\in\LL_{ad}(X,\tau)$.
\end{lm}

\doc The fact that $(RC(X,\tau),\rho,\BBBB)$ is an LCA is
proved in \cite[Example 40]{VDDB2}. The rest can be easily
checked. \sqs

\begin{lm}\label{lead4}
Let $(X,\tau)$ be a Tychonoff space and
$(RC(X),\rho,\BBBB)\in\LL_{ad}(X)$. Then there exists a unique
separated local proximity space $(X,\b,\BB)$ on $(X,\tau)$ such
that $\BBBB=RC(X)\cap\BB$ and $\b_{|RC(X)}=\rho$. In more details,
we set $\BB=\{M\sbe X\st \ex B\in\BBBB$ such that $M\sbe B\}$, and
for every $M,N\sbe X$, we put $ M(-\b)N\iff \fa B\in\BB\ \ex
F,G\in RC(X)$ such that $M\cap B\sbe\int_X(F),\ N\cap
B\sbe\int_X(G)$ and $F(-\rho)G$.
\end{lm}

\doc The proof that $(X,\b,\BB)$ is a separated local
proximity space on $(X,\tau)$ is straightforward; for verifying
the axiom (BC1) we use Theorem \ref{lead5}. The uniqueness follows
from Lemma \ref{lead1}. \sqs




\begin{defi}\label{locsd}
\rm Let $X$ be a locally compact Hausdorff space. We will denote
by $\LL_a(X)$ the set of all LCAs of the form $(RC(X),\rho,\BBBB)$
which satisfy the following conditions:

\smallskip

\noindent(LA1) for every $F,G\in RC(X)$,  $F\cap G\nes$ implies
$F\rho G$;

\smallskip

\noindent(LA2) $CR(X)\sbe\BBBB$;

\smallskip

\noindent(LA3) for every $F\in RC(X)$ and every $G\in CR(X)$,
$F\rho G$ implies $F\cap G\nes$.

\smallskip

If
$(A,\rho_i,\BBBB_i)\in\LL_a(X)$, where $i=1,2$, we set
$(A,\rho_1,\BBBB_1)\preceq_l (A,\rho_2,\BBBB_2)$     if
$\rho_2\sbe\rho_1$ and $\BBBB_2\sbe\BBBB_1$.
\end{defi}

\begin{theorem}\label{locset}
Let $(X,\tau)$ be a locally compact Hausdorff space. Then there
exists an isomorphism $\mu$  between the posets $(\LL(X),\le)$ and
$(\LL_a(X),\preceq_l)$.
\end{theorem}

\doc  It follows immediately from Theorem \ref{lead5}. \sqs

\begin{nota}\label{algcomp}
\rm If $(A,\rho,\BBBB)$ is a CLCA then we will write
$\rho\sbe_{\BBBB} C$ provided that $C$ is a normal contact
relation on $A$ satisfying the following conditions:

\noindent(RC1) $\rho\sbe C$, and

\noindent(RC2) for every $a\in A$ and every $b\in\BBBB$, $aCb$
implies $a\rho b$.

\smallskip

\noindent If $\rho\sbe_{\BBBB} C_1$ and $\rho\sbe_{\BBBB} C_2$
then we will write $C_1\preceq_c C_2$ iff $C_2\sbe C_1$.
\end{nota}

\begin{cor}\label{compset}
Let $(X,\tau)$ be a locally compact Hausdorff space and set
$(A,\rho,\BBBB)=(RC(X),\rho_X,CR(X))$. Then  there exists an
isomorphism  $$\mu_c:(\KK(X),\le)\lra(\KK_a(X),\preceq_c),$$ where
$\KK_a(X)$ is the set of all normal contact relations $C$ on $A$
such that $\rho\sbe_{\BBBB} C$ (see \ref{algcomp} for the
notations).
\end{cor}

\doc It follows immediately from Theorem \ref{locset}.
 \sqs

\begin{pro}\label{alphabeta}
Let $(X,\tau)$ be a locally compact non-compact Hausdorff space
and set $(A,\rho,\BBBB)=(RC(X),\rho_X,CR(X))$. Then $C_\rho$ (see
\ref{Alexprn} for this notation) is the smallest element of the
poset $(\KK_a(X),\preceq_c)$; hence, if $(\a X,\a)$ is the
Alexandroff (one-point) compactification  of $X$ then $\mu_c([(\a
X,\a)])=C_\rho$ (see Corollary \ref{compset} for $\mu_c$).
Further,  the poset $(\KK_a(X),\preceq_c)$ has a greatest element
$C_{\b\rho}$; it is defined as follows: for every $a,b\in A$,
$a(-C_{\b\rho}) b$ iff there exists a set $\{c_d\in A\st
d\in\DDDD\}$  such that:

\noindent(1) $a\llx c_d\llx b^*$, for all $d\in\DDDD$, and

\noindent(2) for any two elements $d_1,d_2$ of $\DDDD$, $d_1<d_2$
implies that $c_{d_1}\llx c_{d_2}$.

\noindent Hence,  if $(\b X,\b)$ is the Stone-\v{C}ech
compactification of $X$  then $\mu_c([(\b X,\b)])=C_{\b\rho}$.
\end{pro}

\doc It is straightforward. \sqs

\begin{rem}\label{rembeta}
\rm The definition of the relation $C_{\b\rho}$ in Proposition
\ref{alphabeta} is given in the language of contact relations. It
is clear that if we use the fact that all happens in a topological
space $X$ then we can  define the relation $C_{\b\rho}$ by setting
for every $a,b\in A$, $a(-C_{\b\rho}) b$ iff $a$ and $b$ are
completely separated.
\end{rem}

\begin{pro}\label{joins}
Let $X$ be a locally compact non-compact Hausdorff space. Set
$(A,\rho,\BBBB)=(RC(X),\rho_X,CR(X))$ and let $\{C_m\st m\in M\}$
be a subset of $\KK_a(X)$ (see \ref{compset} for $\KK_a(X)$). For
every $a,b\in A$, put $a(-C)b$ iff there exists a set $\{c_d\in
A\st d\in\DDDD\}$ such that:

\noindent(1) $a\ll_{C_m} c_d\ll_{C_m} b^*$, for all $d\in\DDDD$
and for each $m\in M$, and

\noindent(2) for any two elements $d_1,d_2$ of $\DDDD$, $d_1<d_2$
implies that $c_{d_1}\ll_{C_m} c_{d_2}$, for every $m\in M$.

\noindent Then $C$ is the supremum in $(\KK_a(X),\preceq_c)$ of
the set $\{C_m\st m\in M\}$.
\end{pro}

\doc The proof is straightforward. \sqs

\section{Extensions over Local Compactifications}

\begin{theorem}\label{zdextcl}
Let, for $i=1,2$, $(X_i,\tau_i)$ be a Tychonoff space, $(Y_i,f_i)$
be a Hausdorff local compactification of $(X_i,\tau_i)$,
$(RC(X_i),\eta_i,\BBBB_i)=\a_{X_i}([(Y_i,f_i)])$ (see (\ref{dw1l})
and (\ref{01l}) for $\a_{X_i}$),
 and $f:X_1\lra X_2$ be a continuous function. Then there exists a continuous function
 $g=L(f):Y_1\lra Y_2$ such that $g\circ f_1=f_2\circ f$ iff $f$ satisfies the
 following conditions:

\smallskip

\noindent{\rm (REQ1)} For every $F,G\in RC(X_2)$,
$\cl_{X_1}(\int_{X_1}(f\inv(F)))\eta_1\cl_{X_1}(\int_{X_1}(f\inv(G)))$
implies that $F\eta_2 G$;

\smallskip

\noindent{\rm (REQ2)} For every $F\in \BBBB_1$ there exists $G\in
\BBBB_2$ such that $f(F)\sbe G$.
\end{theorem}

 \noindent{\em First Proof.}~  Set $\BB_i=\{M\sbe X_i\st \ex
B\in\BBBB_i$ such that $M\sbe B\}$, where $i=1,2$. For every
$M,N\sbe X_i$, set $M(-\eta\ap_i)N\iff \fa B\in\BB_i\ \ex F,G\in
RC(X_i)$ such that $M\cap B\sbe\int_{X_i}(F),\ N\cap
B\sbe\int_{X_i}(G)$ and $F(-\eta_i)G$, where $i=1,2$. Then, by
\ref{lead4}, for $i=1,2$, $(X_i,\eta\ap_i,\BB_i)$ is the unique
separated local proximity space such that $\BB_i\cap
RC(X_i)=\BBBB_i$ and $(\eta\ap_i)_{|RC(X_i)}=\eta_i$. So, if we
prove that $f:(X_1,\eta\ap_1,\BB_1)\lra(X_2,\eta\ap_2,\BB_2)$ is
equicontinuous iff it satisfies conditions (REQ1) and (REQ2), our
assertion will follow from Leader's Theorem \ref{Leader}.

It is easy to see that $f$ satisfies condition (EQ2) iff it
satisfies condition (REQ2). Let $f$ be an equicontinuous function,
$F_1,F_2\in RC(X_2)$ and
$$\cl_{X_1}(\int_{X_1}(f\inv(F_1)))\eta_1\cl_{X_1}(\int_{X_1}(f\inv(F_2))).$$
Then
$\cl_{X_1}(\int_{X_1}(f\inv(F_1)))\eta\ap_1\cl_{X_1}(\int_{X_1}(f\inv(F_2)))$
and thus $$f(\cl_{X_1}(\int_{X_1}(f\inv(F_1))))\eta\ap_2
f(\cl_{X_1}(\int_{X_1}(f\inv(F_2)))).$$ Since, for $i=1,2$,
$f(\cl_{X_1}(\int_{X_1}(f\inv(F_i)))\sbe
\cl_{X_2}f((\int_{X_1}(f\inv(F_i))))\sbe F_i$, we get that
$F_1\eta\ap_2 F_2$ and, therefore, $F_1\eta_2 F_2$. Hence, $f$
satisfies condition (REQ1). So, every equicontinuous function
satisfies conditions (REQ1) and (REQ2). Conversely, let $f$
satisfies conditions (REQ1) and (REQ2), $M,N\sbe X_1$ and
$M\eta\ap_1 N$. Then there exists $B\in\BB_1$ such that for every
$H_1,H_2\in RC(X_1)$ with $M\cap B\sbe\int_{X_1}(H_1)$ and $N\cap
B\sbe\int_{X_1}(H_2)$, $H_1\eta_1 H_2$ holds. Suppose that
$f(M)(-\eta\ap_2)f(N)$. Then, for every $C\in\BB_2$ there exist
$F,G\in RC(X_2)$ such that $f(M)\cap C\sbe \int_{X_2}(F)$,
$f(N)\cap C\sbe \int_{X_2}(G)$ and $F(-\eta_2)G$. Since condition
(REQ2) implies condition (EQ2), we have that $f(B)\in\BB_2$. Thus
there exist $F,G\in RC(X_2)$ such that $f(M)\cap f(B)\sbe
\int_{X_2}(F)$, $f(N)\cap f(B)\sbe \int_{X_2}(G)$ and
$F(-\eta_2)G$. Then $M\cap
B\sbe\int_{X_1}(\cl_{X_1}(\int_{X_1}(f\inv(F))))$ and $N\cap
B\sbe\int_{X_1}(\cl_{X_1}(\int_{X_1}(f\inv(G))))$. Hence, by
(REQ1), $F\eta_2 G$ holds, a contradiction. Therefore,
$f(M)\eta\ap_2 f(N)$. Thus, $f$ is an equicontinuous function.

\smallskip

\noindent{\em Second Proof.}~ In the first proof we used the
Leader Local Compactification Theorem \ref{Leader}. We will now
give another proof which does not use Leader's theorem. Hence, by
the First Proof, it will imply the second part of Leader Theorem
\ref{Leader}. The more important thing is that the method of this
new proof will be used later on for the proof of our Main Theorem
\ref{zdextcmainl}.

\noindent($\Rightarrow$) Let there exists a continuous function
$g:Y_1\lra Y_2$ such that $g\circ f_1=f_2\circ f$. Then, using the
notations
 of (\ref{001l}), we have, by the proof of Theorem
 \ref{lead5}, that the maps
$r_i=r_{(Y_i,f_i)}$ are LCA-isomorphisms, $i=1,2$. Set, for
$i=1,2$, $e_i=(r_i)\inv$ and $\rho_i=\rho_{Y_i}$. Then, by
\ref{isombool}, for every $F\in RC(X_i)$ and $i=1,2$,
$e_i(F)=\cl_{Y_i}(f_i(F))$.  Let $\p_g=\LAM^t(g)$ (see Theorem
\ref{lccont}), i.e.
\begin{equation}\label{0icl}
\ \ \p_g:(RC(Y_2),\rho_2,CR(Y_2))\lra(RC(Y_1),\rho_2,CR(Y_1)), \
G\mapsto \cl_{Y_1}(g\inv(\int_{Y_2}(G))).
\end{equation}
Set also
\begin{equation}\label{0icl1}
\p_f=r_1\circ\p_g\circ
e_2:(RC(X_2),\eta_2,\BBBB_2)\lra(RC(X_1),\eta_1,\BBBB_1).
\end{equation}
 We will prove that
\begin{equation}\label{psifl}
\p_f(G)=\cl_{X_1}(f\inv(\int_{X_2}(G))), \mbox{ for every } G\in
RC(X_2).
\end{equation}
Indeed, let $G\in RC(X_2)$. Then
$\p_f(G)=(f_1)\inv(\cl_{Y_1}(g\inv(\int_{Y_2}(\cl_{Y_2}(f_2(G))))))=
\cl_{X_1}((f_1)\inv(g\inv(\int_{Y_2}(\cl_{Y_2}(f_2(G))))))$. It is
easy to see that
$$(f_2)\inv(\int_{Y_2}(\cl_{Y_2}(f_2(G))))=\int_{X_2}(G).$$ Thus we
obtain that $(f_1)\inv(g\inv(\int_{Y_2}(\cl_{Y_2}(f_2(G)))))=
\{x\in X_1\st (g\circ
f_1)(x)\in\int_{Y_2}(\cl_{Y_2}(f_2(G)))\}=\{x\in X_1\st
f_2(f(x))\in\int_{Y_2}(\cl_{Y_2}(f_2(G)))\} =\{x\in X_1\st f(x)\in
(f_2)\inv(\int_{Y_2}(\cl_{Y_2}(f_2(G))))\}= \{x\in X_1\st f(x)\in
\int_{X_2}(G)\}=f\inv(\int_{X_2}(G))$. Now it becomes clear that
(\ref{psifl}) holds.

Since, by Theorem \ref{lccont}, $\p_g$ is a $\DHLC$-morphism, we
get that $\p_f$ is a $\DHLC$-morphism. Therefore, by (DLC4), for
every $F\in \BBBB_1$ there exists $G\in \BBBB_2$ such that
$F\sbe\p_f(G)$. Since $g$ is continuous, we get that $f$ is
continuous. Thus $f(F)\sbe f(\cl_{X_1}(f\inv(\int_{X_2}(G)))\sbe
\cl_{X_2}(f(f\inv(\int_{X_2}(G))))\sbe G$. Hence, condition (REQ2)
is  checked. Let now $F,G\in RC(X_1)$ and $F\ll_{\eta_1}G$. Then,
by condition (DLC3S), $(\p_f(F^*))^*\ll_{\eta_2}\p_f(G)$. It is
easy to see now that condition (REQ1) is also fulfilled.

\smallskip

\noindent($\Leftarrow$) Let $f$ be a function satisfying
conditions (REQ1) and (REQ2). Set
$\p_f:(RC(X_2),\eta_2,\BBBB_2)\lra(RC(X_1),\eta_1,\BBBB_1)$,
$G\mapsto \cl_{X_1}(f\inv(\int_{X_2}(G)))$. Then it is easy to
check that $\p_f$ is a $\DHLC$-morphism. Put $g=\LAM^a(\p_f)$.
Then $g$ is a continuous function and
$g:\LAM^a(RC(X_1),\eta_1,\BBBB_1)\lra\LAM^a(RC(X_2),\eta_2,\BBBB_2)$,
i.e. $g:Y_1\lra Y_2$  (see the prof of Theorem \ref{lead5}). We
will show that $g\circ f_1=f_2\circ f$. Let $x\in X_1$. Then
$g(f_1(x))= g(\s_x)$ and $f_2(f(x))=\s_{f(x)}$. By Theorem
\ref{lccont}, we have that $g(\s_x)\cap\BBBB_2=\{G\in\BBBB_2\st\fa
F\in RC(X_2), (G\ll_{\eta_2} F)\rightarrow(x\in\p_f(F))\}$.  We
will prove that $\{G\in\BBBB_2\st\fa F\in RC(X_2), (G\ll_{\eta_2}
F)\rightarrow(x\in\p_f(F))\}=\{G\in \BBBB_2\st f(x)\in G\}$. This
will imply, by \ref{bbcl},
 the desired equality. So, let $G\in\BBBB_2$
 and $f(x)\in G$. Let $F\in RC(X_2)$ and $G\ll_{\eta_2} F$. Using condition
(A1), we get that $G\ll_{\rho_{X_2}} F$, i.e. that $G\sbe
\int_{X_2}(F)$. Thus we obtain that $x\in f\inv(G)\sbe
f\inv(\int_{X_2}(F)\sbe\p_f(F)$. Conversely, let $G\in\BBBB_2\cap
g(\s_x)$. Suppose that $f(x)\nin G$. Then $f(x)\in X_2\stm
G=\int_{X_2}(G^*)$. By condition (A2), there exists $F\in \BBBB_2$
such that $f(x)\in\int_{X_2}(F)$ and $F\ll_{\rho_2} G^*$. Then
$G\ll_{\rho_2} F^*$. Hence $x\in\p_f(F^*)$. Since
$f(x)\in\int_{X_2}(F)=X_2\stm F^*$, we get a contradiction.
Therefore, $f(x)\in G$. Thus, $g\circ f_1=f_2\circ f$. \sqs

It is natural to write $f:(X_1,RC(X_1),\rho_1,\BBBB_1)\lra
(X_2,RC(X_2),\rho_2,\BBBB_2)$ when we have a situation like that
which is described in Theorem \ref{zdextcl}. Then, in analogy with
the Leader's equicontinuous functions (see Leader's Theorem
\ref{Leader}), the continuous functions
$f:(X_1,RC(X_1),\rho_1,\BBBB_1)\lra (X_2,RC(X_2),\rho_2,\BBBB_2)$
which satisfy conditions (REQ1) and (REQ2) will be called {\em
R-equicontinuous functions}.

Recall that a function $f:X\lra Y$ is called {\em skeletal}\/
(\cite{MR}) if
\begin{equation}\label{ske}
\int(f\inv(\cl (V)))\sbe\cl(f\inv(V))
\end{equation}
for every open subset  $V$  of $Y$. Recall also the following
three results:

\begin{lm}\label{mrro}{\rm (\cite{Di2})}
Let $f:X\lra Y$ be a continuous map. Then the
following conditions are equivalent:\\
(a) $f$ is a skeletal map;\\
(b) For every $F\in RC(X)$, $\cl(f(F))\in RC(Y)$.
\end{lm}

\begin{lm}\label{skelnewcor}{\rm (\cite{Di6})}
A  continuous map $f:X\lra Y$, where  $X$ and $Y$ are topological
spaces, is skeletal iff for every open dense subset $V$ of $Y$,
$\cl_X(f\inv(V))= X$ holds.
\end{lm}

\begin{lm}\label{skelnewnew}{\rm (\cite{Di6})}
Let $(X_i,\tau_i)$, $i=1,2$, be two topological spaces,
$(Y_i,f_i)$ be some extensions of $(X_i,\tau_i)$, $i=1,2$,
$f:X_1\lra X_2$ and $g:Y_1\lra Y_2$ be two continuous functions
such that $g\circ f_1=f_2\circ f$. Then $g$ is skeletal iff $f$ is
skeletal.
\end{lm}

We are now ready to prove the main result of this paper:

\begin{theorem}\label{zdextcmainl}
Let, for $i=1,2$, $(X_i,\tau_i)$ be a Tychonoff space, $(Y_i,f_i)$
be a Hausdorff local compactification of $(X_i,\tau_i)$,
$(RC(X_i),\eta_i,\BBBB_i)=\a_{X_i}([(Y_i,f_i)])$ (see (\ref{dw1l})
and (\ref{01l}) for $\a_{X_i}$),
  $f:(X_1,RC(X_1),\rho_1,\BBBB_1)\lra
(X_2,RC(X_2),\rho_2,\BBBB_2)$
   be an R-equicontinuous function and  $g=L(f):Y_1\lra Y_2$ be the continuous
   function such that $g\circ f_1=f_2\circ f$ (its existence is guaranteed by
  Theorem \ref{zdextcl}). Then:

  \smallskip

\noindent(a) $g$ is skeletal iff $f$ is skeletal;

  \smallskip

\noindent(b) $g$ is an open map iff $f$  is a skeletal map and
satisfies the following condition:

\smallskip

\noindent{\rm(O)}  $\fa F\in \BBBB_1$ and  $\fa G\in RC(X_1),
(F\ll_{\rho_1} G)\rightarrow
(\cl_{X_2}(f(F))\ll_{\rho_2}\cl_{X_2}(f(G)))$;

\noindent(b${}\ap$) $g$ is an open map iff $f$   satisfies the
following condition:

\smallskip

\noindent{\rm(O1)}  $\fa F\in \BBBB_1$ and  $\fa G\in RC(X_1),
(F\ll_{\rho_1} G)\rightarrow
(\cl_{X_2}(f(F))\ll_{\rho_2\ap}\cl_{X_2}(f(G)))$, where
$\rho_2\ap$ is the unique separated local proximity on $X_2$ such
that $(\rho_2\ap)_{|RC(X_2)}=\rho_2$ (see \ref{lead4});

\noindent(b$''$) $g$ is an open map iff $f$   satisfies the
following condition:

\smallskip

\noindent{\rm(O2)}  $\fa A\sbe X_1$ such that there exists $F\in
\BBBB_1$ with $A\sbe F$, and $\fa B\sbe X_1$, $(A\ll_{\rho_1\ap}
B)\rightarrow (f(A)\ll_{\rho_2\ap}\cl_{X_2}(f(B)))$, where, for
$i=1,2$, $\rho_i\ap$ is the unique separated local proximity on
$X_i$ such that $(\rho_i\ap)_{|RC(X_i)}=\rho_i$ (see \ref{lead4});

\noindent(c) $g$ is a perfect map iff $f$ satisfies the following
condition:

\smallskip

\noindent{\rm(P)} For every $G\in \BBBB_2$,
$\cl_{X_1}(f\inv(\int_{X_2}(G)))\in \BBBB_1$ holds;

\smallskip

\noindent(d) $\cl_{Y_2}(g(Y_1))=Y_2$ iff\/
$\cl_{X_2}(f(X_1))=X_2$;

\smallskip

\noindent(e) $g$ is an injection iff $f$ satisfies the following
condition:

\smallskip

\noindent{\rm(I)} For every $F_1,F_2\in \BBBB_1$ such that
$F_1(-\rho_1) F_2$ there exist $G_1,G_2\in \BBBB_2$ with
$G_1\ll_{\rho_2} G_2$, $F_1\sbe\cl_{X_1}(f\inv(\int_{X_2}(G_2)))$
and $\cl_{X_1}(f\inv(\int_{X_2}(G_2)))(-\rho_1)F_2$;

  \smallskip

\noindent(f) $g$ is an open injection iff  $f$ satisfies condition
{\rm (O1)} (or, equivalently, $f$ is skeletal and satisfies
condition (O)) and the following one:

\smallskip

\noindent{\rm(OI)}  $\fa F\in RC(X_1)\ \ex G\in RC(X_2)$ such that
$F=\cl_{X_1}(f\inv(\int_{X_2}(G)))$;

  \smallskip

\noindent(g) $g$ is a perfect surjection iff $f$ satisfies
condition {\rm (P)} and $\cl_{X_2}(f(X_1))=X_2$.
\end{theorem}

\doc Set $\p_g=\LAM^t(g)$ (see Theorem \ref{lccont}). Then
$\p_g:RC(Y_2)\lra RC(Y_1)$, $G\mapsto
\cl_{Y_1}(g\inv(\int_{Y_2}(G)))$. Set also $\p_f: RC(X_2)\lra
RC(X_1)$, $F\mapsto \cl_{X_1}(f\inv(\int_{X_2}(F)))$.  Then,
(\ref{0icl}), (\ref{0icl1}) and (\ref{psifl}) imply that
$\p_f=r_1\circ\p_g\circ e_2$.

\smallskip

\noindent(a) It follows from Lemma \ref{skelnewnew}.

  \smallskip

\noindent(b) Since every open map is skeletal, we get, using (a),
that if $g$ is an open map then $f$ is skeletal. So, we can
suppose that $f$ is skeletal. Then, as it follows from the proof
of \cite[Theorem 2.11]{Di2}, $\p_f$ is a complete Boolean
homomorphism. Thus, by \cite[(33)]{Di2}, the map $\p_f$ has a left
adjoint $\p^f:RC(X_1)\lra RC(X_2), \ F\mapsto\cl_{X_2}(f(F))$.
Further, by the proof of Theorem
 \ref{lead5}, the maps
$r_1:(RC(Y_1),\rho_{Y_1},CR(Y_1))\lra (RC(X_1),\rho_1,\BBBB_1), \
G\mapsto f_1\inv(G)$ and $e_2=r_2\inv$ are LCA-isomorphisms.
Hence, (\ref{0icl1}) and \ref{maintheorem} imply that $g$ is an
open map iff the map $f$ (is skeletal and) satisfies the following
condition:

\smallskip

\noindent{\rm(O${}\ap$)}  $\fa F\in \BBBB_1$ and  $\fa G\in
RC(X_2), (\p^f(F)\rho_2 G)\rightarrow(F\rho_1\p_f(G))$.

\smallskip

\noindent It is easy to see that condition (O${}\ap$) is
equivalent to the following one:

\smallskip

\noindent{\rm(O$''$)}  $\fa F\in \BBBB_1$ and  $\fa G\in RC(X_2),
(F\ll_{\rho_1}\p_f(G))\rightarrow(\p^f(F)\ll_{\rho_2}G)$.

\smallskip

\noindent  We will prove that condition (O$''$) is equivalent to
condition (O). Indeed, let $F\in\BBBB_1$, $G\in RC(X_1)$ and
$F\ll_{\rho_1}G$. Since $f$ is skeletal, the map $\p^f$ exists.
Set $H=\p^f(G)$. Thus, $H\in RC(X_2)$ and
$\p_f(H)=\p_f(\p^f(G))\spe G$. Therefore $F\ll_{\rho_1}\p_f(H)$.
Then (O$''$) implies that $\p^f(F)\ll_{\rho_2}H$, i.e.
$\p^f(F)\ll_{\rho_2}\p^f(G)$. So, condition (O) is satisfied.

Conversely, let $f$ satisfies condition (O), $F\in \BBBB_1$, $
G\in RC(X_2)$ and $F\ll_{\rho_1}\p_f(G)$. Then, by (O),
$\p^f(F)\ll_{\rho_2}\p^f(\p_f(G))$. Since $\p^f(\p_f(G))\sbe G$,
we get that $\p^f(F)\ll_{\rho_2}G$. Thus, condition (O$''$) is
fulfilled.

\smallskip

\noindent(b${}\ap$) Having in mind Lemma \ref{mrro}, we need only
to show that if $f$ satisfies condition (O1) then $f$ is a
skeletal map. So, let $f$ satisfies condition (O1), $V$ be an open
dense subset of $X_2$ and $G=\cl_{X_1}(f\inv(V))$. Then $G\in
RC(X_1)$. Suppose that $G\neq X_1$. Then there exists $x\in
X_1\stm G$. Clearly, $f_1(x)\nin\cl_{Y_1}(f_1(G))$. Since $Y_1$ is
locally compact and Hausdorff, we get that there exists
$F\in\BBBB_1$ such that $x\in F$ and $F(-\rho_1)G$. Thus
$F\ll_{\rho_1}G^*$. Therefore, by (O1),
$\cl_{X_2}(f(F))\ll_{\rho_2\ap}\cl_{X_2}(f(G^*))$. Set $U=X_1\stm
G$. Since $f$ is continuous, we have that
$H=\cl_{X_2}(f(G^*))=\cl_{X_2}(f(U))\sbe \cl_{X_2}(f(X_1\stm
f\inv(V)))=\cl_{X_2}(f(X_1)\stm V)\sbe X_2\stm V$. Thus
$H\ap=\cl_{X_2}(X_2\stm H)\spe\cl_{X_2}(V)=X_2$. We get that
$\cl_{X_2}(f(F))(-\rho_2\ap) X_2$, a contradiction. Thus,
$f\inv(V)$ is dense in $X_1$. Then Lemma \ref{skelnewcor} implies
that $f$ is a skeletal map.

\smallskip

\noindent(b$''$) It is enough to show that conditions (O1) and
(O2) are equivalent. Set $\BBBB_1\ap=\{A\sbe X_1\st \ex
F\in\BBBB_1$ such that $A\sbe F\}$.

Let $f$ satisfies condition (O1), $A\in\BBBB_1\ap, B\sbe X_1$ and
$A\ll_{\rho_1\ap} B$. Then $A(-\rho_1\ap)(X_1\stm B)$. Thus
$\cl_{Y_1}(f_1(A))\cap\cl_{Y_1}(f_1(X_1\stm B))=\ems$. Since
$A\in\BBBB_1\ap$, we have that $\cl_{Y_1}(f_1(A))$ is a compact
subset of $Y_1$. Using the fact that $Y_1$ is a locally compact
Hausdorff space, we get that there exist $F\in RC(X_1)$ and $U\in
RO(X_1)$ such that $A\sbe F$, $X_1\stm B\sbe U$ and $F(-\rho_1\ap)
U$. Set $G=X_1\stm U$. Then $G\in RC(X_1)$ and $F\ll_{\rho_1}G$.
Thus, by (O1), $\cl_{X_2}(f(F))\ll_{\rho_2\ap}\cl_{X_2}(f(G))$.
Since $G\sbe B$, we get that $f(A)\ll_{\rho_2\ap}\cl_{X_2}(f(B))$.
So, $f$ satisfies condition (O2). Obviously, condition (O2)
implies condition (O1). Therefore, conditions (O1) and (O2) are
equivalent.

\smallskip

\noindent(c) By \cite[Theorem 3.7.18]{E2}, $g$ is a perfect map
iff $\p_g$ satisfies the following condition: for every $G\in
CR(Y_2)$, $\p_g(G)\in CR(Y_1)$ holds. Having in mind the proof of
Theorem \ref{lead5} and (\ref{0icl1}), we get that $g$ is a
perfect map iff $f$ satisfies condition (P).

\smallskip

\noindent(d) This is obvious.

\smallskip

\noindent(e) Using again  (\ref{0icl1}), our assertion follows
from \cite[Theorem 3.16]{Di5}.

\smallskip

\noindent(f) It follows from (b), (\ref{0icl1}), and \cite[Theorem
3.23]{Di5}.

\smallskip

\noindent(g) It follows from (c) and (d).
 \sqs

Recall that a continuous map $f:X\lra Y$ is called {\em
quasi-open\/} (\cite{MP}) if for every non-empty open subset $U$
of $X$, $\int(f(U))\nes$ holds. As it is shown in \cite{Di2}, if
$X$ is regular and Hausdorff, and $f:X\lra Y$ is a closed map,
then $f$ is quasi-open iff $f$ is skeletal. This fact and Theorem
\ref{zdextcmainl} imply the following two corollaries:

\begin{cor}\label{zdextcmainccl}
Let $(X_1,\d_1)$, $(X_2,\d_2)$ be two Efremovi\v{c} proximity
spaces, $(cX_i,c_i)=L(X_i,\d_i)$ (see \ref{Leader} for this
notation) be the Hausdorff compactification of $(X_i,\tau_{\d_i})$
corresponding, by the Smirnov Compactification Theorem \cite{Sm2},
to the Efremovi\v{c} proximity space $(X_i,\d_i)$, where $i=1,2$,
$f:(X_1,\d_1)\lra (X_2,\d_2)$ be a proximally continuous function,
    and $g=L(f):cX_1\lra cX_2$ be the continuous function such
that $g\circ c_1=c_2\circ f$ (see \ref{Leader} for its existence).
Then:

  \smallskip

\noindent(a) $g$ is quasi-open iff $f$ is skeletal;

  \smallskip

\noindent(b){\rm (V. Z. Poljakov \cite{Po})} $g$ is an open map
iff $f$ satisfies the following condition:

\smallskip

\noindent{\rm(OC)} For every $A,B\sbe X_1$ such that
$A\ll_{\d_1}B$, $f(A)\ll_{\d_2}\cl_{X_2}(f(B))$ holds.
\end{cor}

\begin{cor}\label{zdextcmaincbl}
Let $X_1$, $X_2$ be two Tychonoff spaces,
  $f:X_1\lra X_2$ be a continuous function and $\b f:\b X_1\lra\b X_2$
  be the extension of $f$ to the Stone-\v{C}ech compactifications of $X_1$ and $X_2$. Then:

  \smallskip

\noindent(a) $\b f$ is quasi-open iff $f$ is skeletal;

  \smallskip

\noindent(b) $\b f$ is an open map iff $f$ satisfies the following
condition:

\smallskip

\noindent{\rm(OB)} For every $A,B\sbe X_1$ which are completely
separated in $X_1$, $f(A)$ and $X_2\stm\cl_{X_2}(f(X_1\stm B))$
are completely separated in $X_2$;

  \smallskip

\noindent(c){\rm (V. Z. Poljakov \cite{Po})} If $X_1$ and $X_2$
are normal spaces then $\b f$ is open iff for every $A,B\sbe X_1$
such that $\cl_{X_1}(A)\sbe\int_{X_1}(B)$,
$\cl_{X_2}(f(A))\sbe\int_{X_2}(\cl_{X_2}(f(B)))$ holds;

 \smallskip

\noindent(d){\rm (A. D. Ta\u{i}manov \cite{AT})} If $X_1$ and
$X_2$ are normal spaces and $f$ is an open and closed map then $\b
f$ is open.
\end{cor}

\begin{rem}\label{rempol}
\rm  In \cite{Po}, after establishing the general result
\ref{zdextcmainccl}(b), V. Z. Poljakov writes  (in the notations
of Corollary \ref{zdextcmaincbl}) that $\b f$ is open iff for
every two completely separated subsets $A$ and $B$ of $X_1$, the
sets $f(A)$ and $\{y\in X_2\st f\inv(y)\sbe B\}$ are completely
separated in $X_2$. Since $\{y\in X_2\st f\inv(y)\sbe
B\}=f^\sharp(B)=X_2\stm f(X_1\stm B)$, we get that Poljakov's
condition implies condition (OB) and thus it is sufficient for the
openness of $\b f$. It is, however, not necessary. Indeed, let
$f:\mathbb{Q}\lra \b\mathbb{Q}$ be the inclusion map. Then $\b
f:\b\mathbb{Q}\lra\b\mathbb{Q}$ is the identity map and hence it
is an open map. Let $A,B\sbe\mathbb{Q}$ and $A,B$ be completely
separated in $\mathbb{Q}$. Then, by Poljakov's condition, the sets
$f(A)$ and $f^\sharp(B)$ are completely separated in
$\b\mathbb{Q}$, i.e.
$\cl_{\b\mathbb{Q}}(f(A))\cap\cl_{\b\mathbb{Q}}(f^\sharp(B))=\ems$.
Since $f^\sharp(B)=f(B)\cup(\b\mathbb{Q}\stm\mathbb{Q})$, we get
that $\cl_{\b\mathbb{Q}}(f^\sharp(B))=\b\mathbb{Q}$. Thus $f(A)$
and $\b\mathbb{Q}$ are completely separated in $\b\mathbb{Q}$, a
contradiction. Hence, the map $f$ does not satisfy Poljakov's
condition.
\end{rem}

\end{document}